\documentclass[11pt]{article}
\usepackage{amsmath,latexsym,epsf,mathabx}
\begin{document}

\newcommand{\nc}{\newcommand}
\def\PP#1#2#3{{\mathrm{Pres}}^{#1}_{#2}{#3}\setcounter{equation}{0}}
\def\ns{$n$-star}\setcounter{equation}{0}
\def\nt{$n$-tilting}\setcounter{equation}{0}
\def\Ht#1#2#3{{{\mathrm{Hom}}_{#1}({#2},{#3})}\setcounter{equation}{0}}
\def\qp#1{{${(#1)}$-quasi-projective}\setcounter{equation}{0}}
\def\mr#1{{{\mathrm{#1}}}\setcounter{equation}{0}}
\def\mc#1{{{\mathcal{#1}}}\setcounter{equation}{0}}
\def\HD{\mr{Hom}_{\mc{D}}}
\def\AdT{\mr{add}_{\mc{D}}}
\def\Kb{\mc{K}^b(\mr{Proj}R)}
\def\kb{\mc{K}^b(\mc{P}_R)}

\newtheorem{Th}{Theorem}[section]
\newtheorem{Def}[Th]{Definition}
\newtheorem{Lem}[Th]{Lemma}
\newtheorem{Pro}[Th]{Proposition}
\newtheorem{Cor}[Th]{Corollary}
\newtheorem{Rem}[Th]{Remark}
\newtheorem{Exm}[Th]{Example}
\newtheorem{Sc}[Th]{}
\def\Pf#1{{\noindent\bf Proof}.\setcounter{equation}{0}}
\def\>#1{{ $\Rightarrow$ }\setcounter{equation}{0}}
\def\<>#1{{ $\Leftrightarrow$ }\setcounter{equation}{0}}
\def\bskip#1{{ \vskip 20pt }\setcounter{equation}{0}}
\def\sskip#1{{ \vskip 5pt }\setcounter{equation}{0}}
\def\mskip#1{{ \vskip 10pt }\setcounter{equation}{0}}
\def\bg#1{\begin{#1}\setcounter{equation}{0}}
\def\ed#1{\end{#1}\setcounter{equation}{0}}
\def\KET{T^{^F\bot}\setcounter{equation}{0}}
\def\KEC{C^{\bot}\setcounter{equation}{0}}
\def\Db{\mc{D}^b(\mr{mod}R)}


\title{\bf  Wakamatsu-silting complexes {\thanks {Supported by the National Science Foundation of China
(Grant Nos. 11171149) and the National Science Foundation for Distinguished Young Scholars of Jiangsu Province (Grant No.BK2012044)}}}
\smallskip
\author{\small {Jiaqun WEI}\\
\small Institute of Mathematics, School of Mathematics Sciences,\\
\small
Nanjing Normal University \\
\small Nanjing 210023, P.R.China\\ \small Email:
weijiaqun@njnu.edu.cn}
\date{}
\maketitle
\baselineskip 15pt
%
%
\begin{abstract}
\vskip 10pt%
We introduce Wakamatsu-silting complexes (resp.,
Wakamatsu-tilting complexes) as a common generalization of both
silting complexes (resp., tilting complexes) and Wakamatsu-tilting
modules. Characterizations of Wakamatsu-silting complexes are given. In particular, we show that a complex $T$ is Wakamatsu-silting if and only if its dual $DT$ is Wakamatsu-silting. It is conjectured that all compact Wakamatsu-silting
complexes are just silting complexes. We prove that the conjecture lies under
the finitistic dimension conjecture.

\mskip\




\end{abstract}
%
\vskip 30pt

\section{Introduction}

%
%
%
\hskip 18pt
%
%
%
%
%
%
%

Throughout this paper, $R$ always  denotes an artin algebra. We
denote by $\mr{mod}R$ the category of all finitely generated left
$R$-modules and by $\mc{D}^b(\mr{mod}R)$ the bounded derived category of
$\mr{mod}R$.  The homotopy category of bounded complexes of finitely
generated projective modules is denoted by $\kb$.

A very important notion in the derived category is the notion of {\it tilting complexes} [Rk], since they characterize derived equivalences, see [K][Rk]. A slightly general notion is the notion of {\it silting complexes}   [AI][KV][Ws]. Recall that a complex $T\in \mc{D}^b(\mr{mod}R)$ is  silting
provided (1) $T\in \kb$, (2) $T$ is semi-selforthogonal, i.e., $\mr{Hom}_{\mc{D}}(T,T[i])=0$ for all $i>0$, and (3)
$T$ generates $\kb$. In term of this notion, a tilting complex is  just a
selforthogonal silting complex, i.e., a silting complex $T$ such that $\mr{Hom}_{\mc{D}}(T,T[i])=0$ for all $i\neq 0$ .

The notion of silting (tilting) complexes is a far generalization of  tilting
modules. Another far generalization of tilting modules is the
notion of {\it Wakamatsu-tilting} modules [Wk1][GRS]. Recall that an
$R$-module $T$ (in $\mr{mod}R$) is Wakamatsu-tilting [Wk1],
if it satisfies (1) $T$ is selforthogonal and (2) there is a long
exact sequence $0\to R\to^{f_0} T_0\to^{f_1} T_1\to^{f_2}\cdots\
$ such that all $T_i\in\mr{add}_{\mc{D}}T$ and all
$\mr{Im}f_i\in {^{\bot_{k>0}}T}$.

Wakamatsu-tilting modules doesn't  induce derived equivalences in
general. However, they are connected with an equivalence between more general categories than derived categories, namely, repetitive equivalences. Here, we say that two artin algebras $R$ and $S$ are {\it repetitive equivalent} if their repetitive algebras $\hat{R}$ and $\hat{S}$ are stably equivalent. By Happel's result [Hb], for an artin algebra $R$, there is a fully faithful embedding of the bounded derived category $\Db$ into the stable module category of the repetitive algebra $\hat{R}$ and this embedding is an equivalence if and only if the global dimension of $R$ is finite. Following ideas in [Wk2], we proved that a Wakamatsu-tilting $R$-module $T$ such that its Auslander class ${_T\mc{X}}$ is covariantly finite always
induces a repetitive equivalence between $R$ and $\mr{End}_RT$ [We]. Note that such class of Wakamatsu-tilting modules
contain all tilting modules and that all Wakamatsu-tilting modules over an algebra of finite representation type satisfy the condition.

Repetitive equivalences are more general than derived equivalences. In fact, by results in [As][Ch][Rk] etc., if two artin algebras are derived
equivalent, then their repeptitive algebras are derived equivalent, and hence stably
equivalent. It follows that  a tilting complex $T$ over $R$
always induces a repetitive equivalence between $R$ and $\mr{End}T$.

The above facts show that both Wakamatsu-tilting modules and tilting complexes
contribute to some parts of characterizations of repetitive
equivalences. This suggests to study the common generalization of these two kinds of objects. Moreover, it is also needed when we consider the Morita theory for repetitive equivalences.

In this paper, we introduce the notion of Wakamatsu-silting  (Wakamatsu-tilting) complexes, which is certainly a common generalization of both Wakamatsu-tilting modules and silting (tilting) complexes. In term of this new notion, Wakamatsu-tilting modules are just modules which are Wakamatsu-silting. Generally, a Wakamatsu-silting complex is not compact in the derived category of $R$. We conjectured that compact Wakamatsu-silting complexes are just silting complexes. In fact, we prove the conjecture provided that the finitistic dimension conjecture holds for the algebra $R$. It is the case if the injective dimension of $R$ (considered as a left $R$-module) is finite, or $R$ is an Igusa-Todorov algebra [Wig]. We provide interesting characterizations of Wakamatsu-silitng complexes. In particular, we show that the notion of Wakamatsu-silting complexes is self-dual, in sense that a complex $T$ is Wakamatsu-silting if and only if $DT$ is Wakamatsu-silting, where $D$ is the usual duality functor for artin algebras. Assume that $F: \mc{D}^b(\mr{mod}R)\rightleftharpoons  \mc{D}^b(\mr{mod}S) :G$ define an derived equivalence, we obtain that $T\in  \mc{D}^b(\mr{mod}R)$ is Wakamatsu-silting if and only if $F(T)$ is Wakamatsu-silting. The above two results also provides us many examples of Wakamatsu-silting complexes and Wakamatsu-tilting complexes other than Wakamatsu-tilting modules and tilting complexes.

\mskip\

 Complexes in the paper are always cochain ones and subcategories
are always full subcategories in $\mc{D}^b(\mr{mod}R)$
closed under quasi-isomorphisms. Let $\mr{I}$ be a finite interval
of integers, we denote by $\mc{D}^{\mr{I}}(\mr{mod}R)$ the
subcategory of all complexes whose homologies concenter in
$\mr{I}$.
 We
denote by $fg: L\to N$ the composition of two homomorphism $f:L\to
M$ and $g:M\to N$. For simple, we use $\mr{Hom}_{\mc{D}}(-,-)$ instead of $\mr{Hom}_{\mc{D}^b(\mr{mod}R)}(-,-)$.

For basic knowledge on triangulated categories,  derived
categories and the tilting theory, we refer to [{Hb}] and the
Handbook of tilting theory [{AHKb}].


%
\bskip\
\vskip 30pt

\section{Auslander classes and co-Auslander classes in derived categories}
\mskip\
%

\def\HD{\mr{Hom}_{\mc{D}}}

\hskip 18pt

This section is devoted to basic properties  about Auslander
classes and co-Auslander classes in derived categories, which will
be needed for our main results in the next section.

We use notions following from [Ws]. For reader's convenience, we
recall some of them.

Let $\mc{C}$ be a subcategory  containing  $0$. The  subcategory $\mc{C}$ is
extension closed if for any triangle $U\to V\to W\to$ with
$U,W\in\mc{C}$, it holds that $V\in \mc{C}$.  It is resolving
(resp., coresolving) if it is further closed under the functor
$[-1]$ (resp., $[1]$). Note that $\mc{C}$ is resolving (resp.,
coresolving) if and only if, for any triangle $U\to V\to W\to$
(resp.,  $W\to V\to U\to$) in $\mc{T}$ with $W\in\mc{C}$, one has
that `\ $U\in\mc{C}\Leftrightarrow V\in\mc{C}$\ '.

For a complex $L$, we say that $L$ has a
$\mc{C}$-resolution (resp., $\mc{C}$-coresolution) with the length
at most $n$ ($n\ge 0$), denoted by $\mc{C}$-res.dim$(L)\le n$
(resp., $\mc{C}$-cores.dim$(L)\le n$),  if there is a series of
triangles $L_{k+1}\to X_k\to L_k\to$ (resp., $L_k\to X_k\to
L_{k+1}\to$), where $0\le k\le n$, such that $L_0=L$, $L_{n+1}=0$
and each $X_k\in \mc{C}$.

There are the following subcategories associated with the subcategory
$\mc{C}$, where $n\ge 0$.

\sskip\

\hskip 20pt $(\hat{\mc{C}})_n:=\{L\in\mc{D}^b(\mr{mod}R)\ |\
\mc{C}$-res.dim$(L)\le n\}$.

\hskip 20pt $(\check{\mc{C}})_n:=\{L\in\mc{D}^b(\mr{mod}R)\ |\
\mc{C}$-cores.dim$(L)\le n\}$.

\hskip 20pt $\hat{\mc{C}}\ \ \ \ :=\{L\in\mc{D}^b(\mr{mod}R)\ |\ L\in
(\hat{\mc{C}})_n$ for some $n\}$.

\hskip 20pt $\check{\mc{C}}\ \ \ \ :=\{L\in\mc{D}^b(\mr{mod}R)\ |\ L\in
(\check{\mc{C}})_n$ for some $n\}$.

\sskip\

Let $\mr{I}$ be a class of integers. We have the following notions.
\sskip\

\hskip 20pt ${\mc{C}}^{\bot_{k\in \mr{I}}}\ :=\{N\in \mc{D}^b(\mr{mod}R)\ |\
\HD(M,N[i])=0$ for all $M\in \mc{C}$ and all $k\in \mr{I}\}$.

\hskip 20pt $^{\bot_{k\in \mr{I}}}{\mc{C}}\ :=\{N\in \mc{D}^b(\mr{mod}R)\ |\
\HD(N,M[i])=0$ for all $M\in \mc{C}$ and all $k\in \mr{I}\}$.
%
%
%
%
%
%
%
%


\sskip\

%
%
%

It is easy to see that, for fixed integer $m$,  the subcategory
${\mc{C}}^{\bot_{k>m}}$ (resp., $^{\bot_{k>m}}\mc{C}$) is
coresolving (resp., resolving) and closed under direct summands.

Let $M$ be a complex. Denote by $\mr{add}_{\mc{D}}M$ the
subcategory of all direct  summands of finite coproducts of $M$.
$M$ is said to be {\it semi-selforthogonal} (resp., {\it selforthogonal})
provided that $M\in M^{\bot_{k>0}}$ (resp., $M\in M^{\bot_{k\neq
0}}$).

\bskip\

\noindent{\bf{Setup}}: {\it We fix that $M$ is a
semi-selforthogonal complex throughout this section.}

\bskip\

Now we introduce the following subcategories  associated with the
semi-selforthogonal complex $M$. Let $\mr{I}$ be a finite interval
of integers such that $M\in \mc{D}^{\mr{I}}(\mr{mod}R)$.

\bg{verse}

 ${_{M}\mc{X}^{\mr{I}}}\ :=\{X\in M^{\bot_{k>0}}\ |$ there are
 triangles $X_{i+1}\to M_i\to X_i\to$, where $X_0=X$  and each $M_i\in \AdT{M}$,  such that
 all $X_i\in M^{\bot_{k>0}} \cap \mc{D}^{\mr{I}}(\mr{mod}R)$ $\}$.

 ${\mc{X}_M^{\mr{I}}}\ \ :=\{X\in {^{\bot_{k>0}}M}\ |$ there are
 triangles $X_{i}\to M_i\to X_{i+1}\to$, where
$X_0=X$ and each $M_i\in \AdT{M}$,  such that  all $X_i\in
{^{\bot_{k>0}}M} \cap \mc{D}^{\mr{I}}(\mr{mod}R)$ $\}$.

\ed{verse}

It follows from the definition that for  each $X\in
{_{M}\mc{X}^{\mr{I}}}$, there is a triangle $X'\to M_X\to X\to $
such that $M_X\in\AdT{M}$ and $X'\in  {_{M}\mc{X}^{\mr{I}}}$.
Also,  for each $X\in  {\mc{X}_{M}^{\mr{I}}}$, there is a triangle
$X\to M_X\to X'\to $ such that $M_X\in\AdT{M}$ and $X'\in
{\mc{X}_{M}^{\mr{I}}}$.

\bskip\
We set

\bskip\

${_{M}\mc{X}^{\flat}}:=\{X\ |\ X\in {_{M}\mc{X}^{\mr{I}}}$
for some finite interval $\mr{I}$ of integers   such that $M\in
\mc{D}^{\mr{I}}(\mr{mod}R)\}$

\bskip\
\noindent and  call it the {\it {Auslander
class}} (related to $M$). Similarly, we set

\bskip\

${\mc{X}_M^{\flat}}:=\{X\ |\ X\in
{\mc{X}_{M}^{\mr{I}}}$ for some finite interval $\mr{I}$ of
integers  such that $M\in \mc{D}^{\mr{I}}(\mr{mod}R)\}$

\bskip\

\noindent  and  call
it the {\it {co-Auslander class}}. It is clear that
$\widehat{\mr{add}_\mc{D}M}\subseteq {_{M}\mc{X}^{\flat}}$ and
that $\widecheck{\mr{add}_\mc{D}M}\subseteq {\mc{X}_{M}^{\flat}}$
(since $M$ is semi-selforthogonal).

We note that, in case $M$ is a selforthogonal module,
${_{M}\mc{X}^{\mr{[0,0]}}}$ and ${\mc{X}_{M}^{\mr{[0,0]}}}$ are
just the classes $_M\mc{X}$ and $\mc{X}_M$ firstly introduced by
Auslander-Reiten [AR]. The class $_M\mc{X}$ is often called
Auslander class in the literature, see for instance [AF].

The following result gives important properties  of Auslander
classed and co-Auslander classes. The proof is similar to [Ws, Proposition 2.2]. Here
we list for reader's convenience.


%
%
%
\bg{Pro}\label{Ap}

$(1)$ The Auslander class ${_{M}\mc{X}^{\flat}}$  is
coresolving  and closed under direct summands.

$(2)$ The co-Auslander class ${\mc{X}_M^{\flat}}$   is resolving
and closed under direct summands.

\ed{Pro}

\Pf. (1) Note that $0\in \AdT{M}$, so it is easy  to see that the
Auslander class ${_{M}\mc{X}^{\flat}}$ is closed under [1] by the
definition.

 The Auslander class ${_{M}\mc{X}^{\flat}}$ is  also closed under extensions. To see this, let $U\to^f V\to^g W\to $ be a triangle with $U,W\in
{_{M}\mc{X}^{\flat}}$. By the definition of ${_{M}\mc{X}^{\flat}}$, one easily see that there is some common finite interval $\mr{I}$ of integers such that $U,W\in
{_{M}\mc{X}^{\mr{I}}} $. Thus, it is sufficient to show that ${_{M}\mc{X}^{\mr{I}}}$ is closed under extensions.

By the definition of ${_{M}\mc{X}^{\mr{I}}}$, we have triangles
$U_{i+1}\to M'_i\to^{u_i} U_i\to $ and $W_{i+1}\to M''_i\to^{w_i}
W_i\to $ with $U_i, W_i\in {_{M}\mc{X}^{\mr{I}}}$ and
$M'_i,M''_i\in\AdT{M}$ for all $i\ge 0$, where $U_0=U$ and
$W_0=W$. Note that $M^{\bot_{k>0}}\cap \mc{D}^{\mr{I}}(\mr{mod}R)$
is obviously closed under extensions, so $V_0:=V\in
M^{\bot_{k>0}}\cap \mc{D}^{\mr{I}}(\mr{mod}R)$ too. Since $U_0\in
M^{\bot_{k>0}}$, the map $w_0$ can be lifted to a map
$\theta\in\HD(M''_0,V_0)$ through the map $g$. Hence we have the
following triangle commutative diagram for some $V_1$.

\sskip\

 \setlength{\unitlength}{0.09in}
 \begin{picture}(50,19)
%
%
                 \put(18,3.4){\vector(0,-1){2}}
                 \put(27,3.4){\vector(0,-1){2}}
                 \put(35,3.4){\vector(0,-1){2}}

 \put(18,5){\makebox(0,0)[c]{$U_0$}}
                             \put(20,5){\vector(1,0){2}}
                             \put(21,6){\makebox(0,0)[c]{$_f$}}
 \put(27,5){\makebox(0,0)[c]{$V_0$}}
                             \put(31,5){\vector(1,0){2}}
                             \put(32,6){\makebox(0,0)[c]{$_g$}}
 \put(35,5){\makebox(0,0)[c]{$W_0$}}
                             \put(37,5){\vector(1,0){2}}

                 \put(18,9){\vector(0,-1){2}}
                       \put(17,8){\makebox(0,0)[c]{$_{u_0}$}}
                 \put(27,9){\vector(0,-1){2}}
                       \put(25,8){\makebox(0,0)[c]{$(_{{\ \theta\ }}^{{u_0}f})$}}
                 \put(35,9){\vector(0,-1){2}}
                       \put(37,8){\makebox(0,0)[c]{$_{w_0}$}}

 \put(18,11){\makebox(0,0)[c]{$M'_0$}}
                             \put(20,11){\vector(1,0){2}}
                             \put(21,12){\makebox(0,0)[c]{${_{(1, 0)}}$}}
 \put(27,11){\makebox(0,0)[c]{$M'_0\oplus M''_0$}}
                             \put(31,11){\vector(1,0){2}}
                             \put(32,12.5){\makebox(0,0)[c]{$({_{1}^{0}})$}}
 \put(35,11){\makebox(0,0)[c]{$M''_0$}}
                             \put(37,11){\vector(1,0){2}}
                             \put(33,10){\vector(-1,-1){4}}
                             \put(30,8){\makebox(0,0)[c]{$_\theta$}}

                 \put(18,14.5){\vector(0,-1){2}}
                 \put(27,14.5){\vector(0,-1){2}}
                 \put(35,14.5){\vector(0,-1){2}}

 \put(18,16){\makebox(0,0)[c]{$U_1$}}
                             \put(20,16){\vector(1,0){2}}
 \put(27,16){\makebox(0,0)[c]{$V_1$}}
                             \put(31,16){\vector(1,0){2}}
 \put(35,16){\makebox(0,0)[c]{$W_1$}}
                             \put(37,16){\vector(1,0){2}}

%
%

\end{picture}

\noindent Repeating the above process to the triangle $U_1\to
V_1\to W_1\to $, where $U_1,V_1\in {_{M}\mc{X}^{\mr{I}}}$, and so
on, we obtain triangles $V_{i+1}\to M_i\to V_i\to $, for some
$V_i$'s, where $V_0=V$, such that all $M_i=M'_i\oplus
M''_i\in\AdT{M}$ and  all $V_i\in M^{\bot_{i>0}}\cap
\mc{D}^{\mr{I}}(\mr{mod}R)$. It follows that $V\in
{_{M}\mc{X}^{\mr{I}}}$, i.e., ${_{M}\mc{X}^{\mr{I}}}$ is closed
under extensions.

Combining the above, we see that ${_{M}\mc{X}^{\flat}}$ is coresolving.


Finally, we prove that ${_{M}\mc{X}^{\flat}}$ is closed under
direct summands. It  is also sufficient to show that the
subcategory ${_{M}\mc{X}^{\mr{I}}}$ is closed under  direct
summands, for any finite interval $\mr{I}$ of integers. Assume
that $V=U\oplus W\in {_{M}\mc{X}^{\mr{I}}}$. Then $U,W\in
M^{\bot_{i>0}}\cap \mc{D}^{\mr{I}}(\mr{mod}R)$, since both
$M^{\bot_{i>0}}$ and $\mc{D}^{\mr{I}}(\mr{mod}R)$ are closed under
direct summands. By the definition, there is a triangle $V_1\to
M_0\to V\to $ with $M_0\in \AdT{M}$ and $V_1\in
{_{M}\mc{X}^{\mr{I}}}$. Then we have the following triangle
commutative diagram, for some $W_1$.

\mskip\

 \setlength{\unitlength}{0.09in}
 \begin{picture}(50,18)

                 \put(18,3.4){\vector(0,-1){2}}
                 \put(27,3.4){\vector(0,-1){2}}
                 \put(35,3.4){\vector(0,-1){2}}

 \put(18,5){\makebox(0,0)[c]{$U$}}
                             \put(21,5){\vector(1,0){2}}
 \put(27,5){\makebox(0,0)[c]{$V$}}
                             \put(30,5){\vector(1,0){2}}
 \put(35,5){\makebox(0,0)[c]{$W$}}
                             \put(37,5){\vector(1,0){2}}

                 \put(18,9){\vector(0,-1){2}}
                 \put(27,9){\vector(0,-1){2}}
                 \put(35,9){\vector(0,-1){2}}

 \put(18,11){\makebox(0,0)[c]{$W_1$}}
                             \put(21,11){\vector(1,0){2}}
 \put(27,11){\makebox(0,0)[c]{$M_0$}}
                             \put(30,11){\vector(1,0){2}}
 \put(35,11){\makebox(0,0)[c]{$W$}}
                             \put(37,11){\vector(1,0){2}}

                 \put(18,14.5){\vector(0,-1){2}}
                 \put(27,14.5){\vector(0,-1){2}}
                 \put(35,14.5){\vector(0,-1){2}}

 \put(18,16){\makebox(0,0)[c]{$V_1$}}
                              \put(21,16){\vector(1,0){2}}
 \put(27,16){\makebox(0,0)[c]{$V_1$}}
                              \put(30,16){\vector(1,0){2}}
 \put(35,16){\makebox(0,0)[c]{$0$}}
                              \put(37,16){\vector(1,0){2}}

\end{picture}
\sskip\

From the diagram we obtain a triangle $W_1\to M_0\to W\to$ with
$M_0\in \AdT{M}$. Note that there is also a triangle $V_1\to
W_1\to U\to $, from which we can construct a new triangle $V_1\to
W_1\oplus W\to U\oplus W (=V)\to$. It follows that $W_1\oplus W\in
{_{M}\mc{X}^{\mr{I}}}$, since $V_1,V\in {_{M}\mc{X}^{\mr{I}}}$ and
${_{M}\mc{X}^{\mr{I}}}$ is closed under extensions. Now repeating
the process to $W_1$, and so on, we obtain triangles $W_{i+1}\to
M''_i\to W_i\to$ with each $M''_i\in {\AdT{M}}$, where $W_0=W$,
for all $i\ge 0$. Clearly   each $W_i\in  M^{\bot_{i>0}}\cap
\mc{D}^{\mr{I}}(\mr{mod}R)$ by the construction, so we obtain that
$W\in {_{M}\mc{X}^{\mr{I}}}$ by the definition. It follows that
${_{M}\mc{X}^{\mr{I}}}$ is closed under direct summands.

(2) The proof is dual to that for (1). \hfill $\Box$


%
\bskip\

Let $\mc{C}$ be a suncategory. We denote    $\langle \mc{C}
\rangle_+:=\{X\ |\ X=C[n]$ for some $C\in \mc{C}$ and some integer
$n\geq 0 \}$ and   $\langle \mc{C} \rangle_-:=\{X\ |\ X=C[n]$ for
some $C\in \mc{C}$ and some integer $n\leq 0 \}$.
\sskip\
\bg{Th}\label{At}

$(1)$ The subcategory $\langle{_{M}\mc{X}^{\flat}}\rangle_-$  is a
triangulated subcategory closed under direct summands.

$(2)$ The subcategory $\langle{\mc{X}_M^{\flat}}\rangle_+$  is a
triangulated subcategory closed under direct summands.

\ed{Th}

\Pf. We prove (2). Dually, one can obtain the proof of (1).

(2) Clearly, $\langle{\mc{X}_{M}^{\flat}}\rangle_+$   is closed
under [1]. Since ${\mc{X}_{M}^{\flat}}$ is closed under direct
summands,  it  is also seasy to see that
$\langle{\mc{X}_{M}^{\flat}}\rangle_+$ is closed under direct
summands.

We show that $\langle{\mc{X}_M^{\flat}}\rangle_+$   is closed
under [-1]. In fact, for any $N\in
\langle{\mc{X}_M^{\flat}}\rangle_+$, we have that $N=X[n]$ for
some $X\in {\mc{X}_M^{\flat}}$ and some $n\geq 0$. Then
$N[-1]=X[n-1]=(X[-1])[n]$. Since ${\mc{X}_M^{\flat}}$ is closed
under [-1] by Proposition \ref{Ap}, we have that $X[-1]\in
{\mc{X}_M^{\flat}}$. Hence $N[-1]=(X[-1])[n]\in
\langle{\mc{X}_M^{\flat}}\rangle_+$, by the definition.

Finally, we prove that $\langle{\mc{X}_M^{\flat}}\rangle_+$   is
closed under extensions. Assume $N'\to N\to N''\to $ be a triangle
with $N',N''\in \langle{\mc{X}_M^{\flat}}\rangle_+$. Let
$N'=X'[i]$ and $N''=X''[j]$ for some $i,j\ge 0$. Take some $n$
such than $n\ge i,j$, then $i-n<0$ and $j-n<0$. Hence, both
$X'[i-n]$ and $X''[j-n]$ are in ${\mc{X}_M^{\flat}}$, since
${\mc{X}_M^{\flat}}$ is closed under [-1]. Obviously we have a
triangle $N'[-n]\to N[-n]\to N''[-n]\to $. But both $N'=X'[i-n]$
and $N''=X''[j-n]$ are in ${\mc{X}_M^{\flat}}$, so $N[-n]\in
{\mc{X}_M^{\flat}}$, since ${\mc{X}_M^{\flat}}$ is closed under
extensions by Proposition \ref{Ap}. It follows that
$N=(N[-n])[n]\in \langle{\mc{X}_M^{\flat}}\rangle_+$. \hfill
$\Box$

%
%

%
\bskip\

It is easy to see that both $\langle{_M\mc{X}^{\flat}}\rangle_-$
and $\langle{\mc{X}_M^{\flat}}\rangle_+$ contain the smallest
triangulated subcategory containing $M$.

We have another characterization of the triangulated subcategory
$\langle{_{M}\mc{X}^{\flat}}\rangle_-$ (resp.,
$\langle{\mc{X}_M^{\flat}}\rangle_+$ ).

%
%
\sskip\
\bg{Pro}\label{Atc}

$(1)$ The triangulated subcategory
$\langle{_{M}\mc{X}^{\flat}}\rangle_-$  coincides with the
subcategory $\widecheck{_{M}\mc{X}^{\flat}}$.

$(2)$  The triangulated subcategory
$\langle{\mc{X}_{M}^{\flat}}\rangle_+$  coincides with the
subcategory $\widehat{\mc{X}_{M}^{\flat}}$.

\ed{Pro}

\Pf. (1) By the definition of $\widecheck{_{M}\mc{X}^{\flat}}$, it
is easy to see that $\widecheck{_{M}\mc{X}^{\flat}}\subseteq
\langle{_{M}\mc{X}^{\flat}}\rangle_-$, since
$\langle{_{M}\mc{X}^{\flat}}\rangle_-$ is a triangulated
subcategory containing  Auslander class ${_{M}\mc{X}^{\flat}}$. It
remains to show that $\langle{_{M}\mc{X}^{\flat}}\rangle_-
\subseteq \widecheck{_{M}\mc{X}^{\flat}}$. Take any $N\in
\langle{_{M}\mc{X}^{\flat}}\rangle_-$, we have that $N=X[-n]$ for
some $X\in {_{M}\mc{X}^{\flat}}$ and some $n\ge 0$. Since we have
triangles $X[i]\to 0\to X[i+1]$,  where $-n\le i\le -1$, and $0,
X[0]=X\in {_{M}\mc{X}^{\flat}}$, we obtain that $N=X[n]\in
\widecheck{_{M}\mc{X}^{\flat}}$.

(2) Dually. \hfill $\Box$

%
\bskip\

The following well known result is similar as the Schanuel's lemma
in module category, see for instance [Kr, Appendix A].

\bg{Lem}\label{Sch}

Assume that there are triangles $X\to^f M\to^g N\to$ and
$X\to^{f'} M'\to^{g'} N'\to$ such that both
$\mr{Hom}_{\mc{D}}(f,M')$ and  $\mr{Hom}_{\mc{D}}(f',M)$ are epi.
Then $M\oplus X'\simeq M'\oplus X$.

\ed{Lem}
%
\sskip\

Using Lemma \ref{Sch}, we can give the  following
characterizations of Auslander class and co-Auslander class.

\bg{Pro}\label{Ac2}

$(1)$ The subcategory
$M^{\bot_{k>0}}\cap\langle{_{M}\mc{X}^{\flat}}\rangle_-$
coincides with ${_{M}\mc{X}^{\flat}}$.

$(2)$ The subcategory
${^{\bot_{k>0}}M}\cap\langle{\mc{X}_{M}^{\flat}}\rangle_+$
coincides with ${\mc{X}_{M}^{\flat}}$.

\ed{Pro}

\Pf. We prove (2) and leave the reader the proof of (1).

(2) Obviously,  ${\mc{X}_{M}^{\flat}}\subseteq
{^{\bot_{k>0}}M}\cap\langle{\mc{X}_{M}^{\flat}}\rangle_+$.

Now, take any $X\in
{^{\bot_{k>0}}M}\cap\langle{\mc{X}_{M}^{\flat}}\rangle_+$, we have
that $X[-n]\in {\mc{X}_{M}^{\mr{I}}}$, for some $n\ge 0$ and some
suitable interval $\mr{I}$, by the definition. We use induction on
$n$ to show that $X\in {\mc{X}_{M}^{\mr{I}}}$. Then we have also
${\mc{X}_{M}^{\flat}}\supseteq
{^{\bot_{k>0}}M}\cap\langle{\mc{X}_{M}^{\flat}}\rangle_+$ and
hence two subcategories coincide with each other.

If $n=0$, then $X=X[0]\in {\mc{X}_{M}^{\mr{I}}}$ obviously. Now
assume that $n>0$ and that $X[-(n-1)]\in {\mc{X}_{M}^{\mr{I}}}$
implies that $X\in {\mc{X}_{M}^{\mr{I}}}$ (the induction
assumption). Since $X[-n]\in {\mc{X}_{M}^{\mr{I}}}$, there is a
triangle $X[-n]\to M'\to N\to$ such that $M'\in
\mr{add}_{\mc{D}}M$ and $N\in {\mc{X}_{M}^{\mr{I}}}$. Note that
there always triangles $X[-n]\to 0\to X[-(n-1)]\to $. Since
$X[-(n-1)], N\in {^{\bot_{k>0}}M}$ and $M',0\in
\mr{add}_{\mc{D}}M$, we see that these two triangles satisfy the
assumptions in Lemma \ref{Sch}. Hence we have that $M'\oplus
X[-(n-1)]\simeq N$. But ${\mc{X}_{M}^{\mr{I}}}$ is closed under
direct summands by Proposition \ref{Ap}, so that $X[-(n-1)]\in
{\mc{X}_{M}^{\mr{I}}}$. It follows that $X\in
{\mc{X}_{M}^{\mr{I}}}$ by the induction assumption. \hfill $\Box$

%
%
%


\bskip\ \mskip\

\section {Wakamatsu-silting complexes}


\def\HD{\mr{Hom}_{\mc{D}}}
\def\Ad{\mr{Add}_{\mc{D}}}
\def\ad{\mr{add}_{\mc{D}}}
\def\D<{\mc{D}^{\le 0}}
\def\Dn{\mc{D}^{-}}
\def\Db{\mc{D}^b(\mr{mod}R)}
\def\Di#1{\mc{D}^{{#1}}(\mr{mod}R)\setcounter{equation}{0}}

%

\def\KL#1{{^{\bot_{i>0}}(#1)}}
\def\KR#1{{{#1}^{\bot_{i>0}}}}
\def\Kg#1{{{#1}^{\bot_{i\gg 0}}}}
\def\PD{\mr{Pres}_{\D<}^n(\Ad{T})}


\hskip 20pt
\mskip\

%

We note firstly that, in our terms, an $R$-module $T$ is
Wakamatsu-tilting if and only if $T$ is selforthogonal  such that $R\in {\mc{X}_T^{[0,0]}}$.

Now we introduce the notion of Wakamatsu-silting complexes as follows.


\bg{Def}\label{wstd}

A complex $T\in \Db$ is said to be {\bf Wakamatsu-silting}
provided that $T$ is semi-selforthogonal and $R\in
\langle{\mc{X}_T^{\flat}}\rangle_+$. $T$ is called {\bf
Wakamatsu-tilting} if furthermore $T$ is selforthogonal.
\ed{Def}%

From the definition, one immediately obtain that every
Wakamatsu-tilting module is Wakamatsu-silting as a stalk complex.
To see that silting complexes are also Wakamatsu-silting, we can
use the following characterization of Wakamatsu-tilting complexes.

\vskip 15pt

\bg{Pro}\label{wstt}

Let $T$ be a complex in $\Db$. Then the following are equivalent.

$(1)$ $T$ is Wakamatsu-silting.

$(2)$ $T$ is semi-selforthogonal and $\kb\subseteq
\langle{\mc{X}_T^{\flat}}\rangle_+$.

$(3)$ $T$ is semi-selforthogonal and there is a silting complex
$X$ such that $X\in \langle{\mc{X}_T^{\flat}}\rangle_+$.

\ed{Pro}%

\Pf. It is followed from the facts that
$\langle{\mc{X}_T^{\flat}}\rangle_+$ is a triangulated subcategory
by Theorem \ref{At} and that any silting complex is always in
$\kb$ and generates $\kb$, which is the smallest triangulated
subcategory contain $R$. \hfill $\Box$


\vskip 15pt

It follows that any silting complex is Wakamatsu-silting, by
Proposition \ref{wstt} (3).

The above result also suggests the following useful corollary, which shows the
property to be Wakamatsu-silting is a derived invariance.

\bg{Cor}\label{wstde}

Let $F: \Db \to \mc{D}^b(\mr{mod}S)$ be a triangle functor which
defines a derived equivalence. Then $T\in\Db$ is Wakamatsu-silting
if and only if $F(T)$ is Wakamatsu-silting.
\ed{Cor}%

\Pf. Since $F$ is a triangle functor defining a derived
equivalence, we see that $T$ is semi-selforthogonal if and only if
$F(T)$ is semi-selforthogonal and that $R\in
\langle{\mc{X}_T^{\flat}}\rangle_+$ implies $F(R)\in
\langle{\mc{X}_{F(T)}^{\flat}}\rangle_+$. But $F(R)$ is a tilting
complex, we then obtain that $F(T)$ is Wakamatsu-silting by
Proposition \ref{wstt} (3). \hfill $\Box$

\vskip 15pt

The corollary provides us interesting  examples of Wakamatsu-silting
(Wakamatsu-tilting) complexes other than Wakamatsu-tilting modules
and silting complexes. For instance, one take a derived
equivalence $F: \Db \to \mc{D}^b(\mr{mod}S)$ and take a
Wakamatsu-tilting module $T\in \mr{mod}R$ which is not tilting (see for instance [Wk2, Section 3]),
then $F(T)$ is a Wakamatsu-silting (moreover, Wakamatsu-tilting) complex
which, in general, is not a module again.

\vskip 10pt

By the definition of Wakamatsu-silting complexes, it is easy to
see that a  complex $T$ is  Wakamatsu-silting if and only if $T[n]$
is Wakamatsu-silting for some/any integer $n$. Thus, up to shifts,
we may assume that $T\in \Di{[-r,0]}$ for some $r\ge 0$. Note also
that $\Di{[-r,0]}=R^{\bot_{k\not\in [-r,0]}}$.


We have the following useful  characterization of
Wakamatsu-silting complexes.

\bg{Th}\label{wstr}

Assume that $T\in\Di{[-r,0]}$ for some suitable  integer $r$. Then
the following are equivalent.

$(1)$ $T$ is Wakamatsu-silting.

$(2)$ $T$ is semi-selforthogonal and $R\in \mc{X}_T^{[-r,0]}$.
\ed{Th}%

\Pf. $(1)\Rightarrow (2)$ Note that $R\in {^{\bot_{k>0}}T}$ since
$T\in \Di{[-r,0]}$, so, together with the assumption $R\in
\langle{\mc{X}_T^{\flat}}\rangle_+$, we obtain that $R\in
{\mc{X}_T^{\mr{I}}}$ for some interval $\mr{I}$ of integers such
that $T\in \Di{\mr{I}}$, by Proposition \ref{Ac2}. Then, according
to the definition, we have triangles $R_i\to T_i\to R_{i+1}\to$,
where $i\ge 0$ and $R_0:=R$, such that each
$T_i\in\mr{add}_{\mc{D}}T$ and all terms in ${^{\bot_{k>0}}T}\cap
\Di{\mr{I}}$. We will show that all these terms are also in
$\Di{\mr{[-r,0]}}$, i.e., in $R^{\bot_{k\not\in [-r,0]}}$.

In fact, since $R$ and each $T_i\in R^{\bot_{k>0}}$ and
$R^{\bot_{k>0}}$ is coresolving, we easily obtain that each
$R_i\in R^{\bot_{k>0}}$. Now note that, by applying the functor
$\mr{Hom}_{\mc{D}}(R,-)$ to these triangles, we obtain that
$\mr{Hom}_{\mc{D}}(R,R_i[k])\simeq
\mr{Hom}_{\mc{D}}(R,R_{i+1}[k-1])$ for all $k<-r$ and all $i\ge
0$, since $\mr{Hom}_{\mc{D}}(R,T[k])=0$ for all $k<-r$ by
assumptions. Since $\mr{I}$ is a finite interval and each $R_i\in
\Di{\mr{I}}$, there is some $t\le -r$ such that all $R_i\in
R^{\bot_{k<t}}$. Thus,  for any $R_i$, we obtain that

\hskip 20pt $\mr{Hom}_{\mc{D}}(R,R_i[k])\simeq
\mr{Hom}_{\mc{D}}(R,R_{i+1}[k-1])\simeq \cdots\simeq
\mr{Hom}_{\mc{D}}(R,R_{i+1+k-t}[t-1])=0$,

\noindent for all $k<-r$. That is, each $R_i\in R^{\bot_{k<-r}}$.
Hence, we finally obtain all terms in these triangles are in
$R^{\bot_{k\not\in [-r,0]}}=\Di{\mr{[-r,0]}}$. Hence, $R\in
\mc{X}_T^{[-r,0]}$.

 $(2)\Rightarrow (1)$ It is obvious by the definition.
 \hfill $\Box$

\vskip 15pt

In particular, the above result helps us characterize  Wakamatsu-tilting modules
in term of Wakamatsu-silting complexes as follows.

\bg{Pro}\label{wstm}

Assume that $T\in \mr{mod}R$. Then $T$ is  a Wakamatsu-tilting
module if and only if $T$ is Wakamatsu-silting (as a stalk
complex).
\ed{Pro}%

\Pf. We have known that Wakamatsu-tilting modules are
Wakamatsu-silting as stalk complexes. Assume now $T$ is an
$R$-module and is Wakamatsu-silting as a stalk complex. Note that
$\mr{mod}R=\Di{[0,0]}$, so we have that $R\in {\mc{X}_T^{[0,0]}}$,
by the previous theorem. Obviously $T$ is also selforthogonal,
since every semi-selforthgonal module is selforthogonal. Hence we
obtain that $T$ is a Wakamatsu-tilting module by the definition.
\hfill $\Box$

\vskip 15pt

It naturally arise another question:  how to characterize silting
complexes in term of Wakamatsu-silting complexes? We have the
following conjecture. Recall that a complex is compact in $\Db$ if
and only if it is in $\kb$.

\vskip 10pt

\noindent {\bf Conjecture} \ \ Every compact Wakamatsu-silting
complex over an artin algebra is silting.

\vskip 10pt

The above conjecture specifies to Wakamatsu-tilting Conjecture
which asserts that every Wakamatsu-tilting modules of finite
projective dimension is tilting.

We show that the above conjecture is also a direct corollary of
the finitistic dimension conjecture. For this, the following
observation is useful. The proof is left to the reader.


\bg{Lem}\label{fl}

Let $R$ be an artin algebra. Then following are equivalent.

$(1)$ The finitistic dimension conjecture holds for $R$.

$(2)$ There is some fixed integer $n\ge 0$ such that any complex
$X\in\kb\cap\Di{[a,c]}$, where $a\le c$ are any integers, has $0$ as terms at positions not in
$[a-n,c]$ $($up to homotopy equivalences in $\kb$$)$.

\ed{Lem}%


\bg{Th}\label{wstf}

Assume that $R$ satisfies the finitistic dimension conjecture
$($for instance, the injective dimension of $_RR$ is finite or $R$
is an Igusa-Todorov algebra$)$. Then every compact
Wakamatsu-silting complex over $R$ is silting.
\ed{Th}%

\Pf. Up to shifts, we may assume that the compact
Wakamatsu-silting complex $T\in \Di{[-r,0]}$ for some $r\ge 0$.
Then we have that $R\in \mc{X}_T^{[-r,0]}$, by Theorem \ref{wstr}. Hence, there are
triangles $R_i\to T_i\to R_{i+1}\to $, where $i\ge 0$ and
$R_0:=R$, such that $T_i\in \mr{add}_{\mc{D}}T$ and all $R_i\in
{^{\bot_{k>0}}T}\cap\Di{[-r,0]}$. By Lemma \ref{fl}, there is some
fixed integer $n\ge 0$ such that any $X\in\kb\cap\Di{[-r,0]}$ has
0 as terms at positions not in $[-r-n,0]$ (up to homotopy
equivalences in $\kb$). Now take $m=r+n$. By applying the functor
$\mr{Hom}_{\mc{D}}(R_{m+1},-)$ to triangles $R_i\to T_i\to
R_{i+1}\to $, where $0\le i\le m$, we obtain that

\vskip 10pt

\hskip 15pt $\mr{Hom}_{\mc{D}}(R_{m+1},R_m[1])\simeq
\mr{Hom}_{\mc{D}}(R_{m+1},R_{m-1}[2])\simeq \cdots\simeq
\mr{Hom}_{\mc{D}}(R_{m+1},R_0[m+1])$,

\vskip 10pt

\noindent because $R_{m+1}\in {^{\bot_{k>0}}T}$. Note that
$T\in\kb$ since it is compact, so we have all $R_i\in\kb$ for
$i\ge 0$. It follows that $R_{m+1}$ has 0 as terms at positions
not in $[-m,0]$, since $R_{m+1}\in \kb\cap\Di{[-r,0]}$. This
implies that $\mr{Hom}_{\mc{D}}(R_{m+1},R[m+1])=0$. Hence we have
that

\vskip 10pt

 $\mr{Hom}_{\mc{D}}(R_{m+1},R_m[1])\simeq
\mr{Hom}_{\mc{D}}(R_{m+1},R_0[m+1])=
\mr{Hom}_{\mc{D}}(R_{m+1},R[m+1])=0$

\vskip 10pt

\noindent and consequently, the
triangle $R_m\to T_m\to R_{m+1}\to $ splits. Thus, $R_m\in
\mr{add}_{\mc{D}}T$ and $T$ generates $R$. It follows that $T$
generates $\kb$ and $T$ is silting. \hfill $\Box$

%
%
%
%
%
%
%
%
%
%
%
%
%
%
%


\vskip 15pt

Note that an $R$-module $T$ is Wakamatsu-tilting if and only if
its dual module, i.e.,  $DT$ is also Wakamatsu-tilting [Wk2][GRS]. Here
$D$ denotes the usual duality functor for artin algebras. We show that it is still the case for Wakamatsu-silting
complexes.

\bg{Th}\label{Dwst}

Let $T$ be a complex in $\Db$. The following are equivalent.

$(1)$ $T$ is Wakamatsu-silting.

$(2)$ $T$ is semi-selforthogonal and $DR\in
\langle{_T\mc{X}}^{\flat}\rangle_-$.

$(3)$ $DT$ $(in\ \mc{D}^b(\mr{mod}R^o))$ is Wakamatsu-silting.
\ed{Th}%

\Pf. $(1)\Rightarrow (2)$ Up to shifts, we may assume that
$T\in\Di{[-r,0]}$ for some $r\ge 0$. Thus $R\in \mc{X}_T^{[-r,0]}$
and cosequently, there are triangles $R_i\to T_i\to R_{i+1}\to $,
where $i\ge 0$ and $R_0:=R$, such that $T_i\in \mr{add}_{\mc{D}}T$
and all $R_i\in {^{\bot_{k>0}}T}\cap\Di{[-r,0]}$. We claim that
$DR[r]\in {_T\mc{X}}^{[-r,0]}$ and then the statement (2) follows.

Note the fact that $DR[r]\in {T^{\bot_{k>0}}}$ since
$T\in\Di{[-r,0]}$, so we always have triangles $K_{j+1}\to
T'_j\to^{f_j} K_j\to $, where $j\ge 0$ and $K_0:=DR$, such that
each $K_j\in {T^{\bot_{k>0}}}$ (see for instance [Ws, Lemma
3.11]). By applying the functor $\mr{Hom}(R,-)$ to these
triangles, we have that $\mr{Hom}(R,K_j[k-1])\simeq \mr{Hom}(R,
K_{j+1}[k])$ for all $k<-r$ and all $j\ge 0$, since $R\in
{^{\bot_{k<-r}}T}$. In this way, we obtain that

\vskip 10pt

\hskip 15pt ${\mr{Hom}(R,
K_{j+1}[k])\simeq \cdots\simeq
\mr{Hom}(R,K_0[k-j-1])=\mr{Hom}(R,(DR[r])[k-j-1])}$.

\vskip 10pt

\noindent But the latter
is $0$, since $DR[r]\in\Di{[-r,0]}$ and $k-j-1<-r$. Therefore, we
get that $K_j\in R^{\bot_{k<-r}}$ for all $j\ge 0$.

It remains to show that $K_j\in R^{\bot_{k>0}}$ for all $j\ge 0$.
By applying the functor $\mr{Hom}(-,K_j)$ to  triangles $R_i\to
T_i\to R_{i+1}\to $, we obtain that $\mr{Hom}(R_i,K_j[k])\simeq
\mr{Hom}(R_{i+1}, K_j[k+1)$ for all $k>0$ and $i\ge 0$, since
$K_j\in {T^{\bot_{k>0}}}$. Then we get that, for all $k>0$,
$$\mr{Hom}(R, K_{j}[k])=\mr{Hom}(R_0, K_{j}[k])\simeq \mr{Hom}(R_1,
K_{j}[k+1])\simeq\cdots\simeq \mr{Hom}(R_j, K_{j}[k+j]).$$

\noindent Now by applying the functor $\mr{Hom}(R_j,-)$ to
triangles $K_{n+1}\to T'_{n}\to^{f_{n}} K_n\to $, $0\le n\le j-1$,
we have that
$$\mr{Hom}(R_j,K_j[k+j])\simeq \mr{Hom}(R_j,
K_{j-1}[k+j-1])\simeq\cdots\simeq
\mr{Hom}(R_j,K_0[k])=\mr{Hom}(R_j,(DR[r])[k])$$

\noindent for all $k>0$, since $R_j\in {^{\bot_{k>0}}T}$. But it
is easy to see that $\mr{Hom}(R_j,(DR[r])[k])=0$ since $DR[r]\in
\Di{[-r,0]}$ and $k>0$. It follows that, for each $j\ge 0$ and all
$k>0$,
$$\mr{Hom}(R, K_{j}[k])\simeq \mr{Hom}(R_j, K_{j}[k+j])\simeq
\mr{Hom}(R_j,(DR[r])[k])=0.$$
Hence we conclude that  $K_j\in R^{\bot_{k>0}}$ for all $j\ge 0$.

$(2)\Rightarrow (3)$ $DT$ is obviously semi-selforthogonal. Note
that
$\langle{\mc{X}}_{DT}^{\flat}\rangle_+=D(\langle{_T\mc{X}}^{\flat}\rangle_-)$,
so we have that $R=DDR\in \langle{\mc{X}}_{DT}^{\flat}\rangle_+$, for $DR\in \langle{_T\mc{X}}^{\flat}\rangle_-$.
Hence $DT$ is Wakamatsu-silting.

$(3)\Rightarrow (1)$ Note we have proved that the dual complex of
a Wakamatsu-silting complex is also  Wakamatsu-silting, so we
obtain that $T$ is Wakamatsu-silting, since $T=DDT$ and $DT$ is
Wakamatsu-silting. \hfill $\Box$


\vskip 15pt

This result also provides more examples of Wakamatsu-silting complexes. Let us define a complex $T$ to be {\it cosilting} if it is the dual of some silting complex. Then we see that any cosilting complex is Wakamatsu-silting and, in general, it is not silting.


\vskip 15pt

Finally, we ask the following natural question.


\vskip 15pt

\noindent {\bf Question}
Let $T$ be a Wakamatsu-tilting complex over an artin algebra $R$. Does it induce a repetitive equivalence between $R$ and $\mr{End}T$ ?


\vskip 15pt

The following result is easy partial answer to the question.

\bg{Pro}\label{wsre}

Let  $T$ be a Wakamatsu-tilting complex over an artin algebra $R$. Assume that $R$ is derived equivalent to an artin algebra $S$ such that the image of $T$ under the equivalence, say $W$, is $(1)$ a Wakamatsu-tilting module whose related Auslander class is covariantly finite or, $(2)$ a cosilting complex. Then $T$ induces a repetitive equivalence between $R$ and $\mr{End}T$.

\ed{Pro}

\Pf. In both cases, $W$ induces a repetitive equivalence between $S$ and $\mr{End}W$. Since $R$ and $S$ are derived equivalent and $W$ is the image of $T$ under the equivalnce, we obtain that $R$ and $S$ is repetitive equivalent and $\mr{End}T\simeq \mr{End}W$. Hence there is a repetitive equivalence between $R$ and $\mr{End}T$. \hfill $\Box$

%
%
%
%
%
%


%
%
%
%
%
%
%
%
%
%
%
%

\vskip 30pt

\noindent {\bf Reference}
{\small
%

%
   \begin{itemize}
\item[{[AI]}]  T. Aihara and O. Iyama, Silting mutation in triangulated categories, J. Lon. Math. Soc. 85 (2012), no.
3, 633-668.
\item[{[AHKb]}] L. Angeleri-H\"{u}gel, D. Happel and H. Krause (eds),
Handbook of Tilting Theory,   London Math. Soc. Lect. Note
Ser. 332 (2007).
\item[{[As]}] H.	Asashiba, A covering technoque for derived equivalence, J. Algebra, 191 (1997), 382-415.
\item[{[AR]}] M. Auslander and  I. Reiten, Applications of contravariantly
finite subcategories, Adv. Math. 86(1991), 111-152.

\item[{[AF]}]  L. L. Avramov and H. Foxby, Ring homomorphisms and finite gorenstein dimension,  Proc. London
Math. Soc. (3) 75 (1997), no. 2, 241-270.
\item[{[Ch]}]  Q. Chen,  Derived equivalence of repetitive algebras, Adv. Math. (Chinese) 37 (2) (2008), 189-196.
\item[{[GRS]}] E.L Green, I Reiten and {\O}. Solberg,
Dualities on generalized Koszul algebras,
Mem. Amer. Math. Soc., 159 (2002), p. 754.

\item[{[Hb]}]   D. Happel, Triangulated Categories in the Representation
Theory of Finite Dimensional Algebras,  London Math. Soc.
Lect. Note Ser.  119  (1988).
\item[{[K]}]   B. Keller, Deriving DG categories. Ann. Sci. Ecole Norm. Sup. 27 (1994), 63-102.
\item[{[KV]}]  B. Keller and D. Vossieck, Aisles in derived categories, Bull. Soc. Math. Belg. Sér. A 40 (1988), no. 2,
239-253.
\item[{[Kr]}]  H. Krause, Cohomologically cofinite complexes, arXiv:1208.4064 [math.AC].
\item[{[Rk]}]   M. Rickard, Morita theory for derived categories, J. Lond. Math. Soc. 39 (2) (1989), 436-456. %
\item[{[Wig]}] J.  Wei,  Finitistic dimension and Igusa-Todorov algebras, Adv. Math. 222 (2009),
2215-2226.
\item[{[Ws]}] J. Wei, Semi-tilting complexes, Israel J. Math. 194 (2013), 871-893.
\item[{[We]}] J.  Wei, Repetitive equivalences and Wakamatsu-tilting modules, manuscript 2013.

\item[{[Wk1]}] T. Wakamatsu,
On modules with trivial self-extensions,
J. Algebra, 114 (1988),  106-114.
\item[{[Wk2]}] T. Wakamatsu,
Stable equivalence for self-injective algebras and a generalization of tilting modules,
J. Algebra, 134 (1990), 298-325.
\end{itemize}
}

\end{document}